\newtheorem{thm}{Theorem}
\newtheorem{lem}[thm]{Lemma}
\newtheorem{cor}[thm]{Corollary}
\newtheorem{dfn}[thm]{Definition}
\newtheorem{pro}[thm]{Proposition}

\newcommand{\mm}{{\mathfrak M}}
\newcommand{\ma}{{\mathfrak A}}
\newcommand{\mb}{{\mathfrak B}}
\def\F{{\cal F}}
\renewcommand{\L}{{\cal L}}
\newcommand{\I}{{\cal I}}
\newcommand{\Y}{{\cal Y}}

\newcommand{\qcf}[1]{\mbox{Q}^{\mbox{\tiny cf}}_{#1}}
\def\otp{\mbox{otp}}

\def\proof{\noindent {\bf Proof.}\hspace{1mm} }
\def\qed{$\Box$\medbreak }
\def\y{\Y^i}
\def\k{\kappa}
\def\gk{\kappa}

\newcommand{\pl}{\lambda}
\newcommand{\pk}{\kappa}

\newcommand{\po}{\omega}

\newcommand{\se}{\subseteq}

\newcommand{\oL}{{\cal P}}

\newcommand{\la}{\langle}
\newcommand{\ra}{\rangle}

\newcommand{\loo}{{L}_{\omega\omega}}
\newcommand{\lko}{{ L}_{\kappa\omega}}
\newcommand{\lkk}{{ L}_{\kappa\kappa}}
\newcommand{\lkl}{{ L}_{\kappa\lambda}}
\renewcommand{\L}{{\cal L}}
\newcommand{\LS}[1]{\mbox{LS}(#1)}
\documentclass[12pt]{article}
\usepackage{amssymb,latexsym}
\author{
Saharon Shelah
\thanks{Research partially supported by the United States-Israel
        Binational Science Foundation. Publication number [ShVa:726] }\\
Institute of Mathematics\\
Hebrew University\\
Jerusalem, Israel\\
shelah@math.huji.ac.il\\
\and
Jouko V\"a\"an\"anen
\thanks{Research partially supported by
grant 40734 of the Academy of Finland.}\\
Department of Mathematics\\
University of Helsinki\\
Helsinki, Finland\\
jouko.vaananen@helsinki.fi}
\title{A Note on Extensions of Infinitary Logic\footnote{We are
indebted to Lauri Hella,
Tapani Hyttinen and  Kerkko Luosto for useful suggestions.}}
\begin{document}

\maketitle

\begin{abstract}

We 
show that a strong form of the so called
Lindstr\"om's Theorem \cite{Lin69} fails to generalize
to  extensions of \(\lko\) and \(\lkk\):
For weakly compact \(\kappa\) there is no strongest  extension of
\(\lko\) with the 
\((\kappa,\kappa)\)-compactness 
property and the L\"owenheim-Skolem theorem down to
\(\kappa\). With an additional set-theoretic assumption,
 there is no strongest  extension of
\(\lkk\) with the 
\((\kappa,\kappa)\)-compactness property and the 
L\"owenheim-Skolem theorem down to \(<\kappa\).

\end{abstract}
 
%%%%%%%%%%%%%%%%%%%%%%%%

%%%%%%%%%%%%%%%%%%%%%%%%

%\section{Introduction}

By a well-known theorem of Lindstr\"om \cite{Lin69}, first order logic
\(\loo\) is the strongest logic which satisifies the compactness theorem
and the downward L\"owenheim-Skolem theorem. For weakly compact \(\kappa\), 
the infinitary logic \(\lko\) satisfies both the  
\((\kappa,\kappa)\)-compactness property
and the  L\"owenheim-Skolem theorem down to \(\kappa\). 
In \cite{Bar74} Jon Barwise pointed  out that \(\lko\)
is not maximal with respect to these properties,
and asked  what is the strongest logic based on
a weakly compact cardinal \(\kappa\) which still satisfies the
\((\kappa,\kappa)\)-compactness 
property and some other natural conditions suggested by \(\kappa\).
We prove (Corollary 5) that 
for weakly compact \(\kappa\) there is no strongest  extension of
\(\lko\) with the 
\((\kappa,\kappa)\)-compactness property and 
the L\"owenheim-Skolem theorem down to \(\kappa\).
This shows that there is no extension of \(\lko\) which would satisify the
most obvious generalization of Lindstr\"om's Theorem.
A stronger result (Theorem 11)
is proved under an additional assumption.

We use the notation and terminology of \cite[Chapter II]{BarFef85} 
as much as possible.
We will work  
with concrete logics such as { first order logic} \(\loo\),
{ infinitary logic} \(\lkl\) and their extensions \(\loo(\{Q_i:i\in I\})\) and
\(\lkl(\{Q_i:i\in I\})\) by generalized quantifiers. 
Therefore it is not at all critical
which definition of a logic one uses as long as these logics are included and 
some basic closure properties 
are respected. We use \(\L\le\L'\) to denote the sublogic relation.
Let \(\oL\) be a property of logics. A logic \(\L^*\) is {\em 
strongest extension of \(\L\) with \(\oL\)}, 
if 
\begin{enumerate}
\item \(\L\le\L^*\),
\item  \(\L^*\) has property \(\oL\),
\end{enumerate}
 and whenever a logic \(\L'\) 
has the properties 1 and 2, then  \(\L'\le\L^*\). 

Let \(\L\) be a logic.
For infinite cardinals \(\pk\) and \(\pl\), 
\(\L\) is \((\pk,\pl)\){\it-compact}
if for all \(\Phi\se\L\) of power  \(\pk\), 
if each subset of \(\Phi\) of cardinality
\(<\pl\) has a model, then \(\Phi\) has a model. 
\(\L\) is {\it \(\pk\)-compact} if it
is \((\pk,\po)\)-compact. \(\kappa\) is
{\it weakly compact for \(\L\) } if \(\L\) is 
\((\pk,\pk)\)-compact. 
\(\L\) is fully compact if it is \(\kappa\)-compact for
all \(\kappa\). 
\(\L\) has the {\it 
L\"owenheim-Skolem property down to \(\kappa\)},
denoted by \(\LS{\pk}\) 
if every  \(\phi\in\L\)  which has a model, has
a model of cardinality \(\le\pk\).
If every sentence \(\phi\in\L\)  which has a model, has
a model of cardinality \(<\pk\), we say that \(\L\) satisfies 
\(\LS{<\pk}\).
\(\otp(R)\) denotes the order-type of the well-ordering
 \(R\).

\medskip

\begin{thm} \cite{Lin69}
  The logic \(\loo\) is the strongest extension of \(\loo\) with  
\(\aleph_0\)-compactness
and \(\LS{\aleph_0}\). 
\end{thm}

Let \(C\) be a non-trivial class of regular cardinals. Let 
\[\qcf{C}xy\phi(x,y,\vec{z})\iff\phi(\cdot,\cdot ,\vec{z})\mbox{ is a linear
order with cofinality in }C.\] By  \cite{She75},
\(\loo(\qcf{C})\) is always fully compact.
For \(C\) an interval we use
the notation \(\qcf{[\kappa,\lambda)}\) and \(\qcf{[\kappa,\lambda]}\).

\begin{pro}
  \label{1} There is no strongest \(\kappa\)-compact extension of \(\loo\).
In fact:
\begin{enumerate}
\item there are  fully compact logics
\(\L_n\), \(n<\omega\), such that 
\(\L_n\le\L_{n+1}\) for all \(n<\omega\), but
no \(\aleph_0\)-compact
logic can extend each \(\L_n\).

\item There is an \(\aleph_0\)-compact logic
\(\L_1\) and a fully compact logic
\(\L_2\) such that no \(\aleph_0\)-compact
logic can extend both \(\L_1\) and \(\L_2\).

\end{enumerate}

\end{pro}

\proof Let \(\L_n=\loo(\{\qcf{[\aleph_\omega,\infty]}\}\cup
\{\qcf{\aleph_l}:l<n\})\). By \cite{She75}, each \(\L_n\) is
fully compact. Clearly, no \(\aleph_0\)-compact
logic can extend each \(\L_n\).

For the second claim, 
let \(\L_1\) be the logic \(\loo(Q_1)\), where
\(Q_1\) is the quantifier ``there exists uncountable many''
introduced by Mostowski \cite{Mos57}.
This logic is \(\aleph_0\)-compact \cite{Fuh64},  
see \cite[Chapter IV]{BarFef85}
for more recent results. Let \(\L_2\) be the logic \(\loo(Q_B)\),
where \(Q_B\) is the quantifier ``there is a branch'' introduced
by Shelah \cite{She78}. More exactly,
\[Q_BxytuM(x)T(y)(t\le u)\]
if and only if
\(\le_T\) is a partial order of \(T\se M\) 
and there are \(D\),\(\le_D\), \(f\) and \(B\) such that:
  \begin{enumerate}
  \item \(\le_D\) is a total order of \(D\se M\)
  \item \(f:\langle T,\le_T\rangle\rightarrow\langle D,\le_D\rangle\)
is strictly increasing
\item \(\forall s\in D\exists p\in T(f(p)=s)\)
\item \(B\se T\) is totally ordered by \(\le_T\)
\item \(\forall b\in B((p\in T \& p\le_T b)\rightarrow p\in B)\)
\item \(\forall s\in D\exists b\in B(s\le_D f(b)).\)
  \end{enumerate}
The reader is referred to \cite{She78} for a proof of the
full compactness of \(\L_2\).

Suppose there were an \(\aleph_0\)-compact logic \(\L\) containing
both \(\L_1\) and \(\L_2\) as a sublogic. 
It is easy to see that 
the class of countable well-orders 
%\(\langle  A,\le\rangle\) 
can be
expressed as a relativized pseudoelementary class in \(\L\).
This contradicts \(\aleph_0\)-compactness of \(\L\). 
\qed

Lauri Hella pointed out that by elaborating the proof of
claim (2) of the above proposition, we can make \(\L_1\)
fully compact. It was proved in \cite{Wac95} that,
assuming GCH, there is no strongest
extension of \(\L_{\omega\omega}\) which is
\(\aleph_0\)-compact. Our proof of (2) of 
the above proposition is essentially the same as a note,
based on a suggestion of Paolo Lipparini, added after
Theorem 8 of \cite{Wac95}.

\begin{pro} Suppose \(\kappa>\aleph_0\).
  There is no strongest extension of \(L_{\kappa^+\omega}\) with
\(\LS{\kappa}\)
\end{pro}

\proof Let 
\(\L_1=L_{\kappa^+\omega}(\qcf{\aleph_0})\)
and \(\L_2=L_{\kappa^+\omega}(\qcf{[\aleph_1,\kappa]})\).
By using standard arguments with elementary chains
of submodels, it is easy to see that both
\(\L_1\) and \(\L_2\) have \(\LS{\kappa}\), but the
consistent sentence
\[R\mbox{ is a linear order with no last element }\wedge\]
\[\neg\qcf{\aleph_0}xyR(x,y)\wedge
\neg\qcf{[\aleph_1,\kappa]}xyR(x,y)\]
has no models of size \(\le\kappa\). \qed

It was proved in \cite{Wac95} that there is no strongest
extension of \(\L_{\omega\omega}\) with \(\LS{\omega}\).

\begin{lem}
  Suppose \(\kappa\) is weakly compact. Then \(\kappa\)
is weakly compact for 
\(L_{\kappa\omega}(\qcf{\{\aleph_0\}})\)
and for \(\L_{\kappa\omega}(\qcf{[\aleph_1,\kappa]})\).
Moreover, if \(\kappa>\omega\), these logics satisfy \(\LS{\kappa}\).
\end{lem}

\proof The claim concerning \(\LS{\kappa}\) is proved with a 
standard
elementary chain argument. We prove the
weak compactness of \(\L_{\kappa\omega}(\qcf{[\aleph_1,\kappa]})\).
The case of \(L_{\kappa\omega}(\qcf{\{\aleph_0\}})\)
is similar, but easier. For this end, suppose
\(T\) is a set of sentences of
\(\L_{\kappa\omega}(\qcf{[\aleph_1,\kappa]})\)
and \(|T|=\kappa\). We may assume \(T\subseteq \kappa\).
If \(\alpha<\kappa\), then we assume that there is a model \(\mm_\alpha
\models T\cap\alpha\). In view of \(\LS{\kappa}\),
it is not a loss of generality to assume that \(\mm_\alpha=
\langle H(\kappa),R_\alpha\rangle\), where
 \(R_\alpha\subseteq\kappa\times\kappa\).
Let \(R(\alpha,\beta,\gamma)\iff R_\alpha(\beta,\gamma)\).
By weak compactness there is a transitive
\(M\) of cardinality \(\kappa\) such that
\[\langle H(\kappa),\epsilon,T,R\rangle\prec_{\lkk}
\langle M,\epsilon,T^*,R^*\rangle\]
and \(\kappa\in M\).
Let \(\mm=\langle M,S\rangle\), where \(S(x,y)\iff
R^*(\kappa,x,y)\). We claim that \(\mm\models T\).
We need only worry about the cofinality-quantifier.
Cofinalities \(<\kappa\) can be expressed in \(\lkk\),
so they are preserved both ways. Therefore also cofinality
\(\kappa\) is preserved, and no other cofinalities can
occur as the models have cardinality \(\kappa\). \qed

Since the logics \(L_{\kappa\omega}(\qcf{\aleph_0}))\)
and \(\L_{\kappa\omega}(\qcf{[\aleph_1,\kappa]})\)
cannot both be a sublogic of a logic with \(\LS{\kappa}\),
we get from the above lemma:

\begin{cor}
  Suppose \(\kappa>\omega\) 
 is weakly compact. 
Then there is no strongest extension of \(\lko\) for which
\(\kappa\) is weakly compact and which has \(\LS{\kappa}\).
\end{cor}

The logic \(\lko\) actually satisfies the property
\(\LS{<\kappa}\) which is stronger than \(\LS{\kappa}\).
To prove a result like the above corollary for the property
\(\LS{<\kappa}\) we have to work a little harder. 
At the same time we extend the proof to extensions of \(\lkk\).
Here the
cofinality quantifiers \(\qcf{C}\) will not help
as \(\qcf{\{\lambda\}}\) is definable in \(\lkk\) for \(\lambda<\kappa\).
Therefore we use more refined  order-type quantifiers.

\begin{dfn}
  Let \(\lkl(Q)\) denote the formal 
extension of \(\lkl\) by the generalized
quantifier symbol \(Qxy\phi(x,y,\vec{z})\). 
If \(\Y\) is a class of ordinals, 
we get a logic \(\lkl(Q,\Y)\) from \(\lkl(Q)\) by defining
the semantics by
\[\ma\models Qxy\phi(x,y,\vec{c})\iff 
\otp(\{\langle a,b\rangle:\ma\models\phi(a,b,\vec{c})\})\in\Y.\]
If \(\phi\in \lkl(Q,\Y)\) and \(\ma\models\phi\),
we say that \(\ma\models\phi\) holds {\em in the 
\(\Y\)-interpretation}.
\end{dfn}

If \(\ma\) is a model, then \[o(\ma,\Y,\kappa,\lambda)\] is the supremum
of all 
\(\otp(\{\langle a,b\rangle:\ma\models\phi(a,b,\vec{c})\})\)
where \(\phi\in \lkl(\Y)\), \(\vec{c}\in A^{<\lambda}\) and
\(\{\langle a,b\rangle:\ma\models\phi(a,b,\vec{c})\}\) is well-ordered.

\begin{lem}
  \label{trivial}
Suppose \(\kappa\ge\lambda\), \(\phi\in\lkl(Q)\), \(\ma\) is a model,
\(\vec{a}\in A^{<\lambda}\), and \(\Y'\cap o(\ma,\Y,\kappa,\lambda)=\Y\).
Then \(\ma\models\phi(\vec{a})\) in the \(\Y\)-interpretation
if and only if
\(\ma\models\phi(\vec{a})\) in the \(\Y'\)-interpretation.
\end{lem}

\proof This is a straightforward induction of the length
of the formula \(\phi\).\qed

\begin{lem}
  \label{triv-ls} 
  \begin{enumerate}
  \item Suppose \(\kappa>\omega\), \(\phi\in\lkk(Q)\), 
and \(\phi \) has a model \(\ma\) in the \(\Y\)-interpretation. Then 
there is a submodel \(\mb\) of \(\ma\) of cardinality
\(\le 2^{\kappa}\) and \(\Y'\subseteq (2^{\kappa})^+\) such that 
\(\Y'\cap\kappa=\Y\) and \(\mb\models\phi\) in the \(\Y'\)-interpretation.

\item  Suppose \(\kappa=\kappa^{<\kappa}\), 
\(T\se\lkk(Q)\), \(|T|\le\kappa\) 
and \(T \) has a model \(\ma\) in the \(\Y\)-interpretation. Then 
for all \(\xi< \kappa^+\)
there is a submodel \(\mb\) of \(\ma\) of cardinality
\(\le \kappa\) and \(\Y'\subseteq \kappa^+\) such that 
\(\Y\cap\xi=\Y'\cap\xi\)
and \(\mb\models T\) in the \(\Y'\)-interpretation.

  \end{enumerate}
\end{lem}

\proof We may assume \(|A|\ge 2^{\kappa}\).
Let us expand \(\ma\) by 
\begin{enumerate}
\item A well-ordering \(\prec\) the order-type of which
exceed all the order-types of well-orderings definable by
subformulas of \(\phi\) with parameters in \(A\).
\item A new predicate \(P\) which contains those elements
\(d\) of \(A\) for which 
\(\otp(\{\langle a,b\rangle:a\prec b\prec d\})\in\Y\)
\item A prediacte \(F\) which codes an isomorphism from each
 well-ordering, definable by a
subformula of \(\phi\) with parameters in \(A\), onto
an initial segment of \(\prec\).
\end{enumerate}
Let \(\langle \ma,\prec,P,F\rangle\) 
be the expanded structure and \(\langle\mb,\prec^*,P^*,F^*\rangle\) an
\(\lkk\)-elementary substructure of it of cardinality
\(\le 2^{\kappa}\). Let 
\[\Y'=\{\otp(\{\langle a,b\rangle\in B^2:a\prec^* b\prec^* d\}:d\in P^*\}\]
It is easy to see that \(\mb\models\phi\)
in the \(\Y'\)-interpretation. \qed

%DEFINITION OF DIAMOND ON WEAKLY COMPACT

Let \(\pi\) be the canonical well-ordering
of ordered triples of ordinals.
We say that a pair \((\delta_1,Z_1)\), where \(Z_1\subseteq\delta_1\)
{\em codes} a pair  \((\delta_2,Z_2)\), where \(Z_2\subseteq\delta_2\),
if there is a bijection \(f:\delta_2\rightarrow\delta_1\)
such that 
\begin{enumerate}
\item \(\delta_1\) is closed under \(\pi\)
\item \(\pi(0,\alpha,\beta)\in Z_1\iff f(\alpha)<f(\beta)\)
\item \(\pi(1,0,\alpha)\in Z_1\iff f(\alpha)\in Z_2\).
\end{enumerate}

\begin{dfn}
  A cardinal \(\kappa\) satisfies
\(\diamondsuit(\mbox{WC})\) if it is weakly compact and  there is a
 sequence \(\langle A_\alpha : \alpha<\kappa\rangle\) such that
\begin{enumerate}
\item \(A_\alpha\subseteq\alpha\) for \(\alpha<\kappa\).
\item \((\forall A\subseteq \kappa)(\{\lambda<\kappa :
A_\lambda=A\cap\lambda\}\in\I^+)\),
where \(\I\) is the weakly compact ideal on \(\kappa\).
\end{enumerate}
\end{dfn}

\begin{pro}
  \begin{enumerate}
  \item If \(\gk\) is measurable \(> \omega\), then \(\kappa\)
satisfies \(\diamondsuit(\mbox{WC})\).
\item If \(\gk\) is weakly compact \(> \omega\), then
there is a generic extension which preserves all cardinals
and in which  \(\kappa\)
satisfies \(\diamondsuit(\mbox{WC})\).
\item If V=L, then every weakly compact cardinal \(>\omega \) satisfies
\(\diamondsuit(\mbox{WC})\).
\end{enumerate}
\end{pro}

%MAIN RESULT

\begin{thm}\label{main}
  Suppose \(\kappa>\omega\) 
 satisfies \(\diamondsuit(\mbox{WC})\) and \(2^\kappa=
\kappa^+\). Then there is no  strongest extension of
\(\lkk\) for which \(\kappa\) is weakly compact and which
has \(\LS{<\kappa}\).

\end{thm}

\proof We shall construct two sets \(\Y^1,\Y^2\subseteq\kappa^+\) such that
\(\kappa\) is  weakly compact for the logics  \(\lkk(Q,\Y^i)\)
and these logics satisfy
\(\LS{<\kappa}\), but no logic containing both \(\lkk(Q,\Y^1)\)
and \(\lkk(Q,\Y^2)\) satisfies \(\LS{<\kappa}\). The sets \(\Y^i\)
are constructed by induction together with ordinals
\(\xi^i_\alpha<\kappa^+\) such
that:
\[\begin{array}{ll}
  \Y^i=\bigcup_{\alpha<\kappa^+}\Y^i_\alpha&\\
  \Y^i_0=\emptyset&\xi^i_0=0\\
   \Y^i_\alpha=\y_\beta\cap \xi^i_\alpha&\mbox{for }\alpha<\beta\\
  \xi^i_\alpha\le\xi^i_\beta&
   \mbox{for }\alpha<\beta\\
  \Y^i_\nu=\bigcup_{\alpha<\nu}\Y^i_\alpha,&\xi^i_\nu=
\bigcup_{\alpha<\nu}\xi^i_\alpha, \mbox{ for } \nu=\cup\nu\\
\Y^1_\alpha\cap\Y^2_\alpha=\emptyset&\mbox{for }\alpha<\kappa\\
\Y^i_\alpha\subseteq \xi^i_\alpha&\mbox{for }\alpha<\kappa^+\\
\end{array}\]

First we define \(\y_\alpha\) for \(\alpha<\kappa\) in such a way
that \(\lkk(\Y^i)\) will in the end have the property \(\LS{<\kappa}\).

Let \(S_1,S_2\) be a partition of cardinals \(<\k\) into
two stationary sets. Let \(\{\phi^i_\nu:\nu\in S_i\}\)
list all \(\lkk(Q)\)-sentences so that each sentence
is listed as \(\phi^i_\nu\) for stationary many
\(\nu\in S_i\).

Suppose \(\alpha=\lambda+1\) and \(\xi^i_\lambda=\lambda\).
Suppose \(\lambda\in S_i\).
\medskip

\noindent {\bf Case 1.}
Suppose that \((\lambda,A_\lambda)\) codes some pair
\((\xi,Z)\).
In this case we let
\[\y_{\alpha}=\y_\lambda\cup(Z\setminus\lambda),\xi^i_\alpha=\xi\]
\[\Y^{3-i}_{\alpha}=\Y^{3-i}_\lambda.\]
\medskip

\noindent {\bf Case 2.} 
Otherwise we let \(\xi^i_{\alpha}=\lambda\), 
\(\Y^i_\alpha=\Y^i_\lambda\),
\(\Y^{3-1}_\alpha=\Y^{3-i}_\lambda\).
\medskip

Suppose then 
\(\alpha=\lambda+2\), \(\xi^i_\lambda=\lambda\in S_i\)
and we have defined \(\xi^i_{\lambda+1}\) and \(\y_{\lambda+1}\). 

\medskip

\noindent {\bf Case 3.} The sentence \(\phi^i_\lambda\) has 
 a model
in the \(\Y\)-interpretation for some \(\Y\subseteq\k^+\) with
\(\Y\cap\xi^i_{\lambda+1}=\y_{\lambda+1}\). 
By Lemma~\ref{triv-ls} part 2,
\(\phi^i_\lambda\) has a model \(\ma\) of cardinality \(<\kappa\) in the 
 \(\Y\)-interpretation for some \(\Y\subseteq\k\) of
cardinality \(<\kappa\) with
\(\Y\cap\xi^i_{\lambda+1}=\y_{\lambda+1}\). Let \(\mu\) be minimal such that
\(\phi^i_\lambda\in {{\cal L}}_{\mu\mu}(\Y)\). Let \(\xi^i_{\lambda+2}=
o(\ma,\Y,\mu,\mu)\)
and \(\y_{\lambda+2}=\Y\). Let \(\Y^{3-i}_{\lambda+2}=\Y^{3-i}_{\lambda+1}\).
\medskip

\noindent {\bf Case 4.} Otherwise \(\xi^i_{\lambda+1}=\xi^i_{\lambda}\),
\(\Y^i_\alpha=\Y^i_{\lambda+1}\),
\(\Y^{3-1}_\alpha=\Y^{3-i}_{\lambda+1}\).
\medskip

Finally for all other \(\alpha\le\kappa\) we let \(\xi^i_\alpha\) and
\(\y_\alpha\) be defined canonically.

This ends the construction of \(\y_\alpha\) for \(\alpha\le\kappa\).
Note that \(\Y^{1}_\kappa\cap\Y^{2}_\kappa=\emptyset\).
Moreover, if \(\phi^i_\nu\) has a model in the \(\Y\)-interpretation
for some \(\Y\supseteq\y_\kappa\),
then, by construction,  \(\phi^i_\nu\) has a model of cardinality \(<\kappa\)
in the \(\y_\kappa\)-interpretation.

Let \(\Y^i_{\kappa+1}=\Y^i_\kappa\cup\{\kappa\}\)
and \(\xi^i_{\kappa+1}=\kappa+2\).
Next we shall define  \(\y_\alpha\) and \(\xi^i_\alpha\)
for \(\gk+1<\alpha<\gk^+\). For this, 
let \(\la T_\alpha:\kappa<\alpha<\gk^+\ra\)
enumerate all  \(\lkk(Q)\)-theories of cardinality \(\le\gk\)
in a language of cardinality \(\le\gk\)
which satisfy the condition that every subset of 
cardinality \(<\kappa\) has a model in the 
\(\Y^i_\kappa\)-interpretation. Here
we use the assumption \(2^\kappa=\kappa^+\).
We may assume \(T_\alpha\se H(\kappa)\) for all \(\alpha\).

Suppose \(\y_\beta\) and \(\xi^i_\beta\) have been defined for 
\(\beta<\alpha\). If \(\alpha=\cup\alpha\), 
\(\y_\alpha\) and \(\xi^i_\alpha\) are defined canonically.
So assume \(\alpha=\beta+1\). 
Let \(T:H(\kappa)\rightarrow H(\kappa)\) be the function
\(T(a)=T_\beta\cap a\). If \(a\in H(\kappa)\), then
\(T(a)\) has a model 
\(\mb_a\)
in the
\(\y_\kappa\)-interpretation. 
By construction, we may assume \(\mb_a\in H(\kappa)\).
Let \(B:H(\kappa)\rightarrow H(\kappa)\) be the function
\(B(a)=\mb_a\).
%%%
%
Let \(Z\se\kappa\) 
code
\((\xi^i_\beta,\y_\beta)\). 
By \(\diamondsuit(\mbox{WC})\),  
\(W=\{\lambda<\kappa :
A_\lambda=Z\cap\lambda\}\in\I^+\),
where \(\I\) is the weakly compact ideal on \(\kappa\).
Let \(A: \gk\rightarrow H(\gk)\) be the function
\(A(\alpha)=A_\alpha\).
By the definition of \(\I\), there
are a transitive set \(M\) and \(A^*,W^*,Y^*,R^*\) such that
\[\la H(\gk),\epsilon,A,W,\y_\kappa,B,T\ra\prec_{\gk\gk}
            \la M,\epsilon,A^*,W^*,Y^*,B^*,T^*\ra\]
and \(\gk\in W^*\). 
Now \(A^*(\kappa)=Z\) and, by construction, 
\(Y^*\cap\xi^i_\beta=\y_\beta\)
%%%%

It is clear now that \(B(\kappa)\) is a model of \(T_\alpha\)
in the \(Y^*\)-interpretation. By Lemma~\ref{triv-ls} there
is a model \(\mb\) of cardinality
\(\le\kappa\) of \(T_\beta\)
in the \(Y^{**}\)-interpretation for some \(Y^{**}\)
with \(Y^{**}\cap\xi^i_\beta=\y_\beta\). Let 
\(\xi^i_\alpha=o(\mb,Y^{**},\kappa,\kappa)\)
and \(\y_\alpha=Y^{**}\cap\xi^i_\alpha\).

Finally, let \(\y=\bigcup_{\alpha<\kappa^+}\y_\alpha\).
\medskip

\noindent{\bf Claim 1.} \(\lkk(\y)\) satisfies the
\(\LS{<\kappa}\)-property.
\medskip

Suppose \(\phi\) is a sentence of \(\lkk(\y)\) with a model.
Let \(\lambda\in S_i\) such that \(\xi^i_\lambda=\lambda\)
and \(\phi^i_\lambda=\phi\). By the construction of
\(\Y^i_{\lambda+2}\) there is a model of \(\phi\) of
cardinality \(<\kappa\).
\medskip

\noindent{\bf Claim 2.} \(\lkk(\y)\) is weakly \(\kappa\)-compact.
\medskip

Suppose \(T\subseteq \lkk(\y)\) is given and every subset of \(T\)
of cardinality \(<\kappa\) has a model in the \(\y\)-interpretation.
Then \(T=T_\alpha\) for some \(\alpha\).
By construction, every subset of \(T_\alpha\) of cardinality \(<\kappa\) 
has a model in
the \(\y\cap\kappa\)-interpretation. Thus the definition of
\(\y_\alpha\) is made so that \(T_\alpha\) has a model 
\(\mb\) in the \(\Y\)-interpretation for some 
\(\Y\) such that \(\Y\cap o(\mb,\Y,\kappa,\kappa)=\y\cap
o(\mb,\Y,\kappa,\kappa)\). Thus by Lemma~\ref{trivial},
\(\mb\models T_\alpha\) in the \(\y\)-interpretation.
The Claim is proved.

We can now finish the proof of the theorem.
In a logic in which both the quantifier \(Q_{\Y^1}\) and \(Q_{\Y^2}\) 
are definable, we can say that the order-type of a well-ordering
is in \(\Y^1\cap\Y^2\). Thus such a logic cannot satisfy \(\LS{<\kappa}\).
\qed

It is interesting to note that a proof like above would not
be possible for the following stronger
L\"owenheim-Skolem property: A {\em filter-family} is a family
\(\F=(\F(A))_{A\ne\emptyset}\), where
\(\F(A)\) is always a filter on the set \(A\). Luosto \cite{Luo92}
defines the concept of a {\em \((\kappa^+,\omega)\)-neat} filter family.
We will not repeat the definition here, its elements are
invariance under bijections, fineness, \(\kappa^+\)-completeness,
normality and upward relativizability (all defined in \cite{Luo92}).
Suppose \(\L\) is a logic of the form
\(\lkl(\vec{Q})\) for some sequence \(\vec{Q}\) of generalized quantifiers.
We say that \(\L\) has the {\em \(\F,\kappa\)-persistency property},
if for all models \(\ma\) and \(B\in\F(A)\), we have \(\ma\restriction B
\prec\ma\). Luosto proves that if \(\L_1\) and \(\L_2\) both satisfy
the \(\F,\kappa\)-persistency property, then there is \(\L_3\)
such that \(\L_1\le\L_3\), \(\L_2\le\L_3\) and \(\L_3\) satisfies
the  \(\F,\kappa\)-persistency property. Lipparini \cite{Lip87}
proves a similar result for families of limit ultrafilters
related closely to compactness.

Tapani Hyttinen pointed out that the assumption
\(2^\kappa=
\kappa^+\) is not needed in Theorem~\ref{main}, if
\(\kappa\) is assumed to be measurable.


\begin{thebibliography}{99}


\bibitem{Bar74}
 {Jon Barwise},
{Axioms for abstract model theory},
 {Ann. Math. Logic},
{7},
{1974},
 {221--265}.

\bibitem{BarFef85}
 {Model-theoretic logics},
 {Barwise, J. and Feferman, S.},
 {Perspectives in Mathematical Logic},
 {Springer-Verlag},
 {New York},
 {1985},
 {xviii+893}.


\bibitem{Fuh64}
 {Gebhard Fuhrken},
 {Skolem-type normal forms for first-order languages with a
             generalized quantifier},
 {Fund. Math.},
 {54},
 {1964},
 {291--302}.


\bibitem{Lin69}
 {Per Lindstr{\"o}m},
 {On extensions of elementary logic},
 {Theoria},
 {35},
 {1969},
 {1--11}.

\bibitem{Lip87}
Paolo Lipparini,
Limit ultraproducts and abstract logics,
Journal of Symbolic Logic vol. 52 (1987), 437--454.

\bibitem{Luo92}
Kerkko Luosto,
Filters in abstract model theory,
Ph.D. Thesis, University of Helsinki, 1992, 81 pages.



\bibitem{MakShe81}
 {Janos Makowsky and Saharon Shelah},
 {The theorems of {B}eth and {C}raig in abstract model theory.
 {I}{I}. {C}ompact logics},
{Archiv f\"ur Mathematische Logik und Grundlagenforschung},
{21},
 {1981},
 {13--35}.


\bibitem{Mos57}
 {Andrzej Mostowski},
 {On a generalization of quantifiers},
 {Fund. Math.},
{44},
 {1957},
 {12--36}



\bibitem{She75}
 {Saharon Shelah},
 {Generalized quantifiers and compact logic},
 {Trans. Amer. Math. Soc.},
 {204},
 {1975},
 {342--364}.

\bibitem{She78}
 {Saharon Shelah},
 {Models with second-order properties. {I}. {B}oolean algebras
     with no definable automorphisms},
 {Annals of Mathematical Logic},
 {14},
 {1978},
 {1},
{57--72}.

\bibitem{Wac95}
{Marek Wac\l awek},
On ordering of the family of logics
with Skolem-{L\"owenheim} property and
countable compactness property,
in: Quantifiers: Logics, Models and Computation, Vol. 2,
(Micha\l Krynicki, Marcin Mostowski and
Les\l aw Szczerba editors), Kluwer
Academic Publishers, Dordrecht, Boston,
London, 1995, pp.229--236.


\end{thebibliography}
\end{document}